\documentclass[10pt]{amsart}
\input{epsf}
\usepackage{amssymb,latexsym}
\topskip  \numberwithin{equation}{section}
\newtheorem{theorem}{Theorem}[section]
\newtheorem{proposition}[theorem]{Proposition}
\newtheorem{corollary}[theorem]{Corollary}

\newtheorem{remark}[theorem]{Remark}

\begin{document}

\pagenumbering{arabic}
\pagestyle{headings}
\def\sof{\hfill\rule{2mm}{2mm}}
\def\ls{\leq}
\def\gs{\geq}
\def\SS{\mathcal S}
\def\qq{{\bold q}}
\def\txx{{\frac1{2\sqrt{x}}}}
\def\sof{\hfill\rule{2mm}{2mm}}
\def\ls{\leq}
\def\gs{\geq}
\def\SS{\mathcal S}
\def\qq{{\bold q}}
\def\aa{{\overline{1}}}
\def\ab{{\overline{2}}}
\def\an{{\overline{n}}}
\def\as{{\overline{s}}}

\title{ {\sc avoiding $2$-letter signed patterns}}

\author[Toufik Mansour and Julain Wset]
{T. Mansour$^{\rm a}$ and J. West$^{\rm b}$}


\maketitle

\begin{center}
{\small $^A$ LaBRI (UMR 5800), Universit\'e Bordeaux 1,
       351 cours de la Lib\'eration, \\
       33405 Talence Cedex, France}\\[4pt]
{\tt toufik@labri.fr}\\ \ \\

{\small $^B$ Malaspina University-College, 900 Fifth Street,
Nanaimo,\\ BC Canada V9R 5S5}\\[4pt]
{\tt westj@mala.bc.ca}
\end{center}
\section*{Abstract}
Let $B_n$ be the hyperoctahedral group; that is, the set of all
signed permutations on $n$ letters, and let $B_n(T)$ be the set of
all signed permutations in $B_n$ which avoids a set $T$ of signed
patterns. In this paper, we find all the cardinalities of the sets
$B_n(T)$ where $T\subseteq B_2$. This allow us to express these
cardinalities via inverse of binomial coefficients, binomial
coefficients, Catalan numbers, and Fibonacci numbers.
\section{Introduction}

Pattern avoidance proved to be a useful language in a variety of
seemingly unrelated problems, from stack sorting \cite{Kn,Rt,W} to
the theory of Kazhdan-Lusztig polynomials ~\cite{Fb},
singularities of Schubert varieties \cite{LS,SCb}, Chebyshev
polynomials \cite{MV1}(references therein), and Rook polynomials
\cite{MV2}. On the other hand, signed pattern avoidance proved to
be a useful language in combinatorial statistics defined in
type-$B$ noncrossing partitions, enumerative combinatorics,
algebraic combinatorics, geometric combinatorics and singularities
of Schubert varieties; see \cite{Be, BK, Cm, FK, MbRs, Rs, Vr}.

{\bf Restricted permutations.} Let $S_{\{a_1,\dots,a_n\}}$ be the
set of all permutations of the numbers $a_1,\dots,a_n$. For
simplicity let us denote by $S_n$ the set $S_{\{1,2,\dots,n\}}$.
Let $\pi\in S_n$ and $\tau\in S_k$ be two permutations. An {\it
occurrence} of $\tau$ in $\pi$ is a subsequence $1\leq
i_1<i_2<\dots<i_k\leq n$ such that $(\pi_{i_1},\dots,\pi_{i_k})$
is order-isomorphic to $\tau$; in such a context $\tau$ is usually
called a {\it pattern}. We say that $\pi$ {\it avoids} $\tau$, or
is $\tau$-{\it avoiding}, if there is no occurrence of $\tau$ in
$\pi$. The set of all $\tau$-avoiding permutations in $S_n$ is
denoted $S_n(\tau)$. For an arbitrary finite collection of
patterns $T$, we say that $\pi$ avoids $T$ if $\pi$ avoids any
$\tau\in T$; the corresponding subset of $S_n$ is denoted
$S_n(T)$.

{\bf Restricted signed permutations .} We will view the elements
of the hyperoctahedral group $B_n$ as signed permutations written
as $\alpha=\alpha_1\alpha_2\dots\alpha_n$ in which each of the
symbols $1,2,\dots,n$ appears, possibly barred. Clearly, the
cardinality of $B_n$ is $2^nn!$. The barring operation is define
by change the symbol $\alpha_i$ to $\overline{\alpha_i}$ and
$\overline{\alpha_i}$ to $\alpha_i$, so it is an involution, and
the absolute value notation means $|\alpha_i|$ is $\alpha_i$ if
the symbol $\alpha_i$ is not barred, otherwise
$\overline{\alpha_i}$.

Now let $\tau\in B_k$, and $\alpha\in B_n$; we say {\em $\alpha$
contains signed pattern $\tau$} or {$\tau$-containing signed
pattern}, if there is sequence of $k$ indices, $1\leq
i_1<i_2<\dots<i_k\leq n$ such that two conditions hold: $(1)$
$\alpha$ with all bars removed contains the pattern $\tau$ with
all bars removed, i.e., $|\alpha_{i_p}|>|\alpha_{i_q}|$ if and
only if $|\tau_p|>|\tau_q|$ for all $k\leq p>q\leq 1$; and $(2)$
$\alpha_{i_j}$ barred if and only if $\tau_j$ barred for all
$1\leq j\leq k$. For example,
$\alpha=21\overline{3}\overline{4}\in B_4$ contains the signed
patterns $\aa\ab$ and $21$. If $\alpha$ not contains signed
pattern $\tau$, then we say $\alpha$ {\em avoids signed pattern}
$\tau$ or $\tau$-{\em avoiding signed pattern}. The set of
$\tau$-avoiding signed permutations in $B_n$ we denote by
$B_n(\tau)$. More generally we define $B_n(T)=\cap_{\tau\in T}
B_n(\tau)$. The cardinality of $B_n(T)$ we denote by $b_n(T)$.

\begin{proposition} {\rm (see \cite[Section~3]{Rs})}
\label{sym}
Let us define, the reversal (i.e., reading the permutation right-to-left:
$\alpha_1\alpha_2\cdots\alpha_n\mapsto\alpha_n\alpha_{n-1}\cdots\alpha_1$),  the barring (i.e.,
$\alpha_1\alpha_2\cdots\alpha_n\mapsto\overline{\alpha_1}\overline{\alpha_2}\cdots\overline{\alpha_n}$)
and the complement (i.e., $\alpha_1\alpha_2\cdots\alpha_n\mapsto\beta_1\beta_2\cdots\beta_n$
where $\beta_i=n+1-\alpha_i$ if $\alpha_i$ not barred, otherwise
$\beta_i=\overline{n+1-|\alpha_i|}$ for all $i$) on $B_n$. Let us denote by $G_b$ the
group which generated by these three symmetric operations. Then every element $g\in G_b$
give a bijection which show if $T$ and $T'$ are both subsets of
signed patterns in $B_n$ such that $T'=g(T)=\{g(\alpha)| \alpha\in T\}$,
then $b_n(T)=b_n(T')$.
\end{proposition}

%
%
In the symmetric group $S_n$, for every $2$-letters pattern $\tau$
the number of $\tau$-avoiding permutations is one, and for every
pattern $\tau\in S_3$ the number of $\tau$-avoiding permutations
is given by the Catalan number ~\cite{Kn}. Simion
~\cite[Section~3]{Rs} proved there are similar results for the
Hyperoctahedral group $B_n$ (generalized by Mansour \cite{M}), for
every $2$-letters signed pattern $\tau$ the number of
$\tau$-avoiding signed permutations is given by $\sum_{j=0}^n
\binom{n}{j}^2j!$. In the present note, similarly as \cite{SS}, we
find all the cardinalities $b_n(T)$ where $T\subseteq B_2$.

The paper is organized as follows. In section $2$ we treat the
case $|T|=1,2$. In sections $3$, $4$, and $5$ we present all
the values $b_n(T)$ where $T\subseteq B_2$ such that $|T|=3$, $|T|=4$,
and $|T|=5,6,7,8$; respectively.
\section{Two signed patterns}
By taking advantage of Proposition \ref{sym}, the question of
determining the values $b_n(\tau)$ for the $8$ choices of one
$2$-letter signed pattern, reduces to $2$ cases, which are
$\tau=12$ and $\tau=1\ab$. Simion \cite[Proposition~3.2]{Rs}
proved for any $n\geq 0$
\begin{equation}
b_n(12)=b_n(1\ab)=\sum_{k=0}^n{n\choose k}^2k!. \label{eq1}
\end{equation}

Additionally, the second question of determining the values
$b_n(\tau,\tau')$ for $28$ choices of two $2$-letters signed
patterns, reduces $8$ cases.

\begin{remark}
In \cite[Proposition~3.4]{Rs} proved $b_n(12,21)=2n!$ and
$b_n(1\ab,\aa2)=(n+1)!$. On the other hand, $b_2(12,21)=6$ and
$b_3(1\ab,\aa2)=22$. Here, we present the correction of these
cases.
\end{remark}

\begin{theorem}
Given $n$ and two $2$-letters signed patterns $\tau$, $\tau'$. The value
$b_n(\tau,\tau')$ satisfies one of the following relations, according to
which orbit (under reversal, barring, complementation) contains the pair
$\tau$, $\tau'$:
\begin{equation}
b_n(\{12,21\})=b_n(\{12,1\ab\})=b_n(\{2\aa,1\ab\})=b_n(\{2\aa,\aa2\})=(n+1)!;
\label{eq2}
\end{equation}
\begin{equation}
b_n(\{12,\aa\ab\})=b_n(\{12,\ab\aa\})={{2n}\choose n}; \label{eq3}
\end{equation}
\begin{equation}
b_n(\{12,\ab1\})=n!+n!\sum\limits_{i=1}^n\left(\dfrac{1}{i}\sum\limits_{j=0}^{i-1}\dfrac{1}{j!}\right);
\label{eq4}
\end{equation}
\begin{equation}
b_n(\{1\ab,\aa2\})=2\sum\limits_{l=1}^n \sum\limits_{\
i_1+i_2+\dots+i_l=n,\ i_j\geq 1\ } \prod\limits_{j=1}^l i_j!.
\label{eq5}
\end{equation}
\end{theorem}
\begin{proof}
$b_n(\{12,21\})=(n+1)!$ yields immediately by use
\cite[Theoren~1]{M}, the other results in \ref{eq2} are holding by
use \cite[Proposition~3.4]{Rs}. By \cite[Example~4.8]{M} we get
\ref{eq3}.

To verify \ref{eq4}, let us consider the number permutations
$\alpha\in B_n(12,\ab1)$. If $\alpha_n$ barred then $\alpha$
avoids $12,\ab1$ if and only if $(\alpha_1,\dots,\alpha_{n-1})$
avoids $12,\ab1$, which means there are $nb_{n-1}(12,\ab1)$ signed
permutations. If $\alpha_n=i$ unbarred then the smaller symbols
$1,\dots,i-1$ must be barred in $\alpha$ and the larger symbols
$i+1,\dots,n$ must be unbarred, hence the smaller can be permuted
and placed in any positions $1,2,\dots,n-1$. This gives
$b_n(12,\ab1)=nb_{n-1}(12,\ab1)+\sum_{i=0}^{n-1} {{n-1}\choose
i}i!$ for $n\geq 1$. Besides $b_0(T)=1$, hence $(4)$ holds by the
principle of induction on $n$.

To verify \ref{eq5}, let us consider $\alpha\in B_n(1\ab,\aa2)$.
So, by induction on $n$ it is easy to prove there exist a
partition $\alpha=(\alpha^1,\dots,\alpha^l)$ such hold the
following conditions:
\begin{enumerate}
\item Every absolute symbol in $\alpha^j$ is greater than every absolute element
        in $\alpha^{j+1}$ for all $j=1,2,\dots l-1$;
\item All the symbols of $\alpha^j$ are either barred or unbarred;
\item The symbols in $\alpha^j$ are barred if and only if the symbols in
    $\alpha^{j+1}$ are unbarred.
\end{enumerate}
Hence, the rest is easy to check.
\end{proof}
\section{Three signed patterns}
By taking advantage of Proposition \ref{sym}, the question of determining
the values $b_n(T)$ where $T\subset B_2$ and $|T|=3$, for the $56$ choices
of three $2$-letters signed patterns, reduces to $10$ cases as follows.
$$\begin{array}{llll}
    T_1=\{12,1\ab,\aa2\};&\ \   T_2=\{12,1\ab,\aa\ab\};&\ \     T_3=\{12,1\ab,21\};&\ \   T_4=\{12,1\ab,2\aa\};\\
    T_5=\{12,1\ab,\ab1\};&\ \   T_6=\{12,1\ab,\ab\aa\};&\ \     T_7=\{12,\aa\ab,21\};&\ \ T_8=\{12,\aa\ab,2\aa\};\\
                 &\ \   T_9=\{12,2\aa,\ab1\};  &\ \     T_{10}=\{1\ab,\aa2,2\aa\}.&
\end{array}$$

\begin{theorem}
Given $n$ and set $T$ of $2$-letters signed patterns such that
$|T|=3$. The value $b_n(T)$ satisfies one of the following
relations, according to which orbit (under reversal, barring,
complementation) contains $T$:
\begin{equation}
b_n(T_1)=\sum\limits_{d=0}^n \sum\limits_{\ i_0+i_1+\dots+i_d=n-d\
} \prod\limits_{j=0}^d i_d!; \label{eq6}
\end{equation}
\begin{equation}
b_n(T_2)=C_{n+1}; \label{eq7}
\end{equation}
\begin{equation}
b_n(T_3)=n!+ n! \sum\limits_{j=1}^{n} \dfrac{1}{j}; \label{eq8}
\end{equation}
\begin{equation}
b_n(T_4)=b_n(T_5)=n!\sum\limits_{j=0}^n \frac{1}{j!}; \label{eq9}
\end{equation}
\begin{equation}
b_n(T_6)=F_{2n+1}; \label{eq10}
\end{equation}
\begin{equation}
b_n(T_7)=n^2+1; \label{eq11}
\end{equation}
\begin{equation}
b_n(T_8)=2^{n+1}-(n+1); \label{eq12}
\end{equation}
\begin{equation}
b_n(T_9)=n!+\sum\limits_{j=1}^n \sum_{p+q=n-j}p!q!,\label{eq13}
\end{equation}
\begin{equation}
b_n(T_9)=n!\sum\limits_{j=0}^n \binom{n}{j}^{-1},\label{eq13a}
\end{equation}
where $C_m$ and $F_m$ are the $m$th Catalan and Fibonacci numbers;
respectively.
\end{theorem}
\subsubsection*{\bf Proof \ref{eq6}} Let $\alpha\in B_n(T_1)$, and let $m_0$ the first
symbol unbarred by reading $\alpha$ from left-to-right. Since
$\alpha$ avoids $12$ we get all the symbols $\alpha_j<m_0$ are
unbarred, and since $\alpha\in B_n(T_1)$ we get that
$\alpha=(\alpha^0,m_0,\beta)$ where $\beta\in B_{m_0-1}(T_1)$ and
$\alpha^0$ is permutation of the symbols
$\overline{m_0+1},\dots,\an$. By the principle of induction, for
any $\alpha\in B_n(T_1)$ there exist $0\leq d\leq n$ and
$\alpha^{j}$ permutations of the symbols
$\overline{m_{j}+1},\overline{m_j+2},\dots,\overline{m_{j-1}-1}$
for all $0\leq j\leq d+1$ where $m_{-1}=n+1$, $m_{d+1}=0$, $0\leq
m_d<m_{d-1}<\dots<m_0\leq n$ such that
$\alpha=(\alpha^0,m_0,\alpha^1,m_1,\dots,\alpha^d,m_d,\alpha^{d+1})$.
The rest easy to check.\qed

\begin{corollary}
\label{extu2} $b_n(T_1\cup\{\aa\ab\})=2^n$ for any $n\geq 0$, and
$b_n(T_1\cup\{\ab\aa\})=1+{{n+1}\choose 2}$ for any $n\geq 2$.
\end{corollary}
\begin{proof}
Immediately by use the argument proof of \ref{eq6}.
\end{proof}
\subsubsection*{\bf Proof \ref{eq7}.}

A {\em split permutation} is a permutation $\pi=(\pi',\pi'')\in
S_n$, where $\pi'$ and $\pi''$ are nonempty such that every entry
of $\pi'$ is greater than every entry of $\pi''$. For example,
$231$, $312$, and $321$ all the split permutations in $S_3$.

We first check that the number of non-splitting $123$-avoiding
permutations in $S_n$ (denoted by $N_n$) is the $(n-1)$th Catalan
number; $C_{n-1}$. This can doubtless be done directly, but we
will offer a somewhat indirect proof by induction. It is easy to
Check a base case. Now, we suppose this property holds for $j<n$,
that is, $N_j=C_{j-1}$ for $j<n$.

Take the $C_n$ $123$-avoiding permutation (see \cite{Kn}) and
class them according to the first (upper-left-most) place where
they split. Since each permutation in $S_n$ thus decomposes into a
direct sum of a non-splitting $123$-avoiding permutation and an
arbitrary $123$-avoiding permutation, we have that
$$C_n=\sum_{j=1}^n N_jC_{n-j}=N_n+\sum_{j=0}^{n-2} C_jC_{n-1-j}.$$
Since we know that the standard catalan recurrence is
$C_n=\sum_{j=0}^{n-1} C_jC_{n-1-j}$, it follows that $N_n=C_{n-1}$
which completing the induction step.

Now, suppose we have a signed permutations $\pi$ avoids $T_2$.
Then the permutation $|\pi|$ must avoids $123$. For considering
any increasing sequence of length $3$, two of the elements must
receive the same sign.

So it remains to consider all $123$-avoiding permutations and
assign signs to their elements, respecting the condition that if
any element is above and to the right of another element, then it
must be coloured $*$ while the one below and to the left must be
coloured $o$.

Take a $123$-avoiding permutation $\pi$, and suppose that one of
its elements can coloured freely, either $*$ or $o$. Then this
element cannot be below and to the left of any other element;
neither can it be above and to the right of any other element.
This means that the remaining elements in the permutation $\pi$
must be located in the other two blocks; i.e. if out
freely-colourable element is denoted by $s$, we have
$\pi=(\pi',s,\pi'')$ where $\pi'$ and $\pi''$ coloured by $o$ such
that each element of $\pi''$ greater than each element of $\pi'$.
This means that the elements in the upper left block must be
exactly sufficient to fill in all the rows above $s$ and all the
columns to the left of $s$. Therefore the number of these rows is
the same as the number of columns, and the $s$ is in fact situated
on the diagonal.  This means that the permutation splits at $s$.
It follows that any permutation which is non-splitting is also
uniquely colourable.

There is one exception, which is that the single 1-permutation can
be coloured in two ways, by $*$ or by $o$.

Now we are ready to find $b_n(T_2)$ as follows. By induction,
$b_j(T_2)=C_{j+1}$ for $j<n$.

We count $b_n(T_2)$ according to the position of the first split.
let $a_n$ be the number of non-splitting signed permutations which
avoid $T_2$; so, $a_1= 2$, while $a_j = N_j = C_{j-1}$ for $j>1$.

Note that if such a signed permutation splits, then each half can
be coloured independently according to the pattern-avoidance
conditions. This means that

$$b_n(T_2) =\sum_{j=1}^{n-1} a_jb_{n-j}+a_n
= 2C_n + \sum_{j=1}^{n-2} C_jC_{n-j} + C_{n-1}=\sum_{j=0}^{n}
C_jC_{n-j}=C_{n+1}.$$\qed
\subsubsection*{\bf Proof \ref{eq8}.}
Let $\alpha\in B_n(T_3)$; if the symbol $n$ unbarred in $\alpha$,
then since $\alpha$ avoids $12$ and $21$ we get that all symbols
of $\alpha$ barred, so there are $n!$ signed permutations. Let
$\alpha_j=\an$; since $\alpha$ avoids $1\ab$, then the symbol
$\alpha_i$ barred for $i\leq j$, so in this case there are
$\sum_{j=1}^n {{n-1}\choose {j-1}}(j-1)! b_{n-j}(T_3)$ signed
permutations. Hence
    $$b_n(T_3)=n!+(n-1)!\sum_{j=1}^{n} \frac{b_{n-j}(T_3)}{(n-j)!}.$$
Let $b'_n=b_n(T_3)/n!$; so $b'_n=1+\frac{1}{n}\sum_{j=0}^{n-1}
b'_j$, which means that $b'_n-b'_{n-1}=\frac{1}{n}$. Besides
$b'_1=2$, hence $b'_n=1+\sum_{j=1}^n \frac{1}{j}$, as claimed in
\ref{eq8}.

\begin{corollary}
\label{extu3}
$b_n(T_3,\aa\ab)=1+{{n+1}\choose 2}$ for $n\geq 0$.
\end{corollary}
\begin{proof}
Immediately by use the argument proof of \ref{eq8}.
\end{proof}
\subsubsection*{\bf Proof \ref{eq9}.}
To verify \ref{eq9} we find $b_n(T_4)$ and $b_n(T_5)$ by two
steps.

{\bf 1.} Let $\alpha\in B_n(T_4)$; since $\alpha$ avoids $12$ and
$1\ab$ we get $\alpha_1=n$ or $\alpha_1=\overline{i}$. In the
first case, since $\alpha$ also avoids $2\aa$ must be
$\alpha=(n,n-1,\dots,1)$. In the second case there are
$b_{n-1}(T_4)$ signed permutations. Therefore,
$b_n(T_4)=1+nb_{n-1}(T_4)$ for all $n\geq3$ with $b_2(T_4)=5$.
Hence by the principle of induction on $n$ we get the formula for
$b_n(T_4)$.
\begin{corollary}
\label{extu4} $b_n(T_4\cup\{\aa\ab\})=b_n(T_4\cup\{\ab\aa\})=2^n$
and $b_n(T_4\cup\{21\})=2\cdot n!$ for all $n\geq1$.
\end{corollary}
\begin{proof}
By the above argument (proof of $b_n(T_4)$ formula) we obtain the
following.
\begin{enumerate}
\item $b_n(T_4\cup\{21\})=nb_{n-1}(T_4\cup\{21\})$ for all $n\geq2$
    with $b_1(T_4\cup\{21\})=2$, hence $b_n(T_4\cup\{21\})=2\cdot n!$
    for $n\geq 1$.
\item   Let $\alpha$ avoids $T_4$ and $\{\aa\ab\}$; similarly there two cases. In the first
    case $\alpha=(n,\dots,2,1)$. In the second case, since $\alpha$ avoids $T_4\cup\{\aa\ab\}$
    we have that $\alpha=(\overline{i},\beta,n,\dots,i+1,\gamma)$ where all symbols of
    $\beta$ barred and decreasing, and all symbols of $\gamma$ unbarred and decreasing.
    So $b_n(T_4\cup\{\aa\ab\})=1+\sum_{i=1}^n 2^{i-1}$ which means
    $b_n(T_4\cup\{\aa\ab\})=2^n$.
\item   Similarly as the second case we get that
    $b_n(T_4\cup\{\ab\aa\})=2^n$.
\end{enumerate}
\end{proof}
{\bf 2.} Let $\alpha\in B_n(T_5)$; similarly as first step we
obtain either $\alpha_1=n$, or $\alpha_1=\overline{i}$. In the
first case there are $b_{n-1}(T_5)$ signed permutations. In the
second case, since $\alpha$ avoids $\ab1$, all the symbols
$1,2,\dots,i-1$ are barred, so there are $\sum_{i=1}^n
{{n-1}\choose{i-1}}(i-1)!b_{n-i}(T_5)$ signed permutations. Hence
$b_n(T_5)=b_{n-1}(T_5)+(n-1)!\sum_{i=0}^n \frac{b_i(T_5)}{i!}$ for
$n\geq 1$. Let $b'_n=b_n(T_5)/n!$, so
$n(b'_n-b'_{n-1})=b'_{n-1}-b'_{n-2}$ for all $n\geq 2$. Besides
$b'_1=2$ and $b'_0=1$, hence $b'_n=\sum_{j=0}^n \frac{1}{j!}$, as
claimed in the second part of \ref{eq9}.
\begin{corollary}
\label{extu5} $b_n(T_5\cup\{\ab\aa\})=2^n$ for all $n\geq 0$.
\end{corollary}
\begin{proof}
By use the above argument (proof on $b_n(T_5)$ formula) we get
that
$$b_n(T_5\cup\{\ab\aa\})=b_{n-1}(T_5\cup\{\ab\aa\})+b_{n-1}(T_5\cup\{\ab\aa\}).$$
Besides, $b_1(T_5\cup\{\ab\aa\})=2$, hence the corollary holds.
\end{proof}
\subsubsection*{\bf Proof \ref{eq10}.}
Let $\alpha\in B_n(T_6)$; it is easy to see that $\alpha_1=n$ or
$\alpha_1=\overline{i}$. In the first case there are
$b_{n-1}(T_6)$ signed permutations. In the second case, since
$\alpha$ avoids $\ab\aa$ we get that, all the symbols $1,\dots
i-1$ are unbarred, since $\alpha$ avoids $12$ then $\alpha$
contains $(i-1,i-2,\dots,1)$, and since $\alpha$ avoids $12$ and
$1\ab$ we get $\alpha=(\overline{i},\beta,i-1,\dots,1)$, hence
there are $b_{n-i}(T_6)$ signed permutations for $1\leq i\leq n$.
Therefore, for all $n\geq 3$
    $$b_n(T_6)=b_{n-1}(T_6)+b_{n-1}(T_6)+b_{n-2}(T_6)+\dots+b_0(T_6),$$
which means that $b_n(T_6)=3b_{n-1}(T_6)-b_{n-2}(T_6)$. Besides
$b_n(T_6)=5$, hence by the principle of induction \ref{eq10}
holds.
\subsubsection*{\bf Proof \ref{eq11}.} Immediately by
\cite[Theoren~4.4]{M}.
\subsubsection*{\bf Proof \ref{eq12}.}
Let $\alpha\in B_n(T_8)$; if $\alpha_1$ unbarred, then it is east
to see that $(\alpha_2,\dots,\alpha_n)$ is partition to two
decreasing subsequences such that, all the symbols
$|\alpha_1|+1,\dots,n$ are barred and the others symbols are
unbarred, hence there are $2^{n-1}$ signed permutations. If
$|\alpha_1|<n$ and $\alpha_1$ barred, then (similarly) there are
$\sum_{i=1}^{n-1} 2^{i-1}$ signed permutations. Finally, if
$\alpha_1=\an$ then by definitions there are $b_{n-1}(T_8)$ signed
permutations, so $b_n(T_8)=b_{n-1}(T_8)+2\cdot 2^{n-1}-1$. Besides
$b_2(T_8)=5$, $b_1(T_8)=2$, and $b_0(T_8)=1$, hence \ref{eq12}
holds.

\begin{corollary}
\label{extu1}
Let $\tau\in\{21,\ab1\}$; $b_n(T_8\cup\{\tau\})=2n$ for all $n\geq 1$.
\end{corollary}
\begin{proof}
Immediately by use the argument proof of \ref{eq12}.
\end{proof}
\subsubsection*{\bf Proof \ref{eq13}.}
To verify \ref{eq13}, let $\alpha\in B_n(T_9)$, $\alpha_i$ the
first entry unbarred in $\alpha$ ($i$ minimal), and let $\alpha_j$
the last entry unbarred in $\alpha$ ($j$ maximal). Since $\alpha$
avoids $12$ all the symbols which are not barred are decreasing,
and since $\alpha$ avoids $2\aa$ and $\ab1$, we can may write
$\alpha=(\beta,\gamma,\delta)$ where $\beta$ is permutations of
the numbers $\aa,\dots,\overline{i}$,
$\gamma=(n-j,n-1-j,\dots,i+1)$, and $\beta$ is permutations of the
numbers $\overline{n-j+1},\dots,\an$. Hence the rest easy to
check.

\subsubsection*{\bf Proof \ref{eq13a}.}
To verify \ref{eq13a}, let $\alpha\in B_n(T_{10})$, and let $j$
the maximal such that $\alpha_j$ is barred. Since $\alpha$ avoids
$T_{10}$, so we can may write $\alpha=(\beta,\gamma)$ where all
symbols of $\beta$ are barred, all symbols of $\gamma$ are
unbarred, and $|\beta_j|>|\gamma_i|$ for all $i,j$. Hence
$b_n(T_{10})=\sum_{j=0}^n (n-j)!j!$, as claimed in the second part
of \ref{eq13a}.
\section{Four signed patterns}
By taking advantage of Proposition \ref{sym}, the question of
determining the values $b_n(T)$ where $T\subset B_2$ and $|T|=4$,
for the $70$ choices of four $2$-letters signed patterns, reduces
to $16$ cases as follows.
$$\begin{array}{lll}
U_1=\{12,1\ab,\aa2,\aa\ab\};&\ \    U_2=\{12,1\ab,\aa2,21\};&\ \        U_3=\{12,1\ab,\aa2,2\aa\};\\
U_4=\{12,1\ab,\aa2,\ab\aa\};&\ \    U_5=\{12,1\ab,\aa\ab,21\};&\ \      U_6=\{12,1\ab,\aa\ab,2\aa\};\\
U_7=\{12,1\ab,\aa\ab,\ab1\};&\ \    U_8=\{12,1\ab,21,2\aa\};&\ \        U_9=\{12,1\ab,21,\ab1\};\\
U_{10}=\{12,1\ab,21,\ab\aa\};&\ \   U_{11}=\{12,1\ab,2\aa,\ab1\};&\ \   U_{12}=\{12,1\ab,2\aa,\ab\aa\};\\
U_{13}=\{12,1\ab,\ab1,\ab\aa\};&\ \ U_{14}=\{12,\aa\ab,21,\ab\aa\};&\ \ U_{15}=\{12,\aa\ab,2\aa,\ab1\};\\
U_{16}=\{1\ab,\aa2,2\aa,\ab1\}.& &
\end{array}$$

\begin{theorem}
\label{thh4} Given $n$ and set $T$ of $2$-letters signed patterns
such that $|T|=4$. The value $b_n(T)$ for $n\geq3$ satisfies one
of the following relations, according to which orbit (under
reversal, barring, complementation) contains $T$:

\noindent $\begin{array}{ll}
(4.1)  &b_n(U_{14})=0;\\
(4.2)\ &b_n(U_{10})=b_n(U_{15})=2n;\\
(4.3)  &b_n(U_4)=b_n(U_5)=1+{{n+1}\choose 2};\\
(4.4)  &b_n(U_1)=b_n(U_6)=b_n(U_7)=b_n(U_{12})=b_n(U_{13})=2^n;\\
(4.5)  &b_n(U_8)=b_n(U_9)=b_n(U_{16})=2n!;\\
(4.6)  &b_n(U_3)=b_n(U_{11})=\sum_{j=0}^n j!;\\
(4.7)  &b_n(U_2)=n!\left( 1+\sum_{j=0}^{n-1}j!(n-1-j)!\right).
\end{array}$
\end{theorem}
\begin{proof}
Immediately by definition $(4.1)$ holds. $(4.2)$ holds by
Corollary \ref{extu1}, and $(4.4)$ holds by two Corollaries
\ref{extu4} and \ref{extu5}.

To verify $(4.3)$, by Corollaries \ref{extu2} and \ref{extu3} it
is sufficient to prove $b_n(U_7)=2^n$. Let $\alpha\in B_n(U_7)$;
if $\alpha_1$ unbarred, then since $\alpha$ avoids $12,1\ab$ we
have $\alpha_1=n$, and in this case there are $b_{n-1}(U_7)$
signed permutations. If $\alpha_1$ barred, then since $\alpha$
avoids $\aa\ab$ all the symbols $|\alpha_1|+1,\dots,n$ are
unbarred and decreasing (since $\alpha$ avoids $12$), and since
$\alpha$ avoids $\ab1$ all the symbols $|\alpha_1|-1,\dots,1$ are
barred and decreasing (since $\alpha$ avoids $\aa\ab$). So there
are $2^{n-1}$ signed permutations. Thus,
$b_n(U_7)=b_{n-1}(U_7)+2^{n-1}$ for $n\geq 1$. Besides
$b_0(U_7)=1$ and $b_1(U_7)=2$, hence $(4.3)$ holds.

To find $b_n(U_9)$ note that $\alpha\in B_n(U_9)$ has two cases, the first is when
    the symbol $n$ unbarred in $\alpha$, since $\alpha$ avoids $12,21$
    all symbols of $\alpha$ are barred, so there are $n!$ signed permutations.
    The second case when the symbol $n$ barred in $\alpha$, since $\alpha$
    avoids $\ab1,1\ab$ all other symbols are barred, so there $n!$ signed permutations.
    Therefore $b_n(U_9)=2n!$ for $n\geq 1$. Similarly $b_n(U_{16})=2n!$ for $n\geq 1$.
    By Corollary \ref{extu4} $b_n(U_8)=2n!$ for $n\geq 1$, hence $(4.5)$ holds.

To find $b_n(U_3)$ note that $\alpha\in B_n(U_3)$ has two cases. The first
    if $\alpha_n$ unbarred, then since $\alpha$ avoids $12,\aa2$ must be $\alpha_n=1$,
    so there are $b_{n-1}(U_3)$ signed permutations. The second case when $\alpha_n$
    is barred, since $\alpha$ avoids $1\ab,2\aa$ we have all symbols of $\alpha$ are
    barred, which means there are $n!$ signed permutations. Accordingly
    $b_n(U_3)=b_{n-1}(U_3)+n!$ for $n\geq 1$, and $b_0(U_3)=1$, hence
    $b_n(U_3)=\sum_{j=0}^n j!$.

To find $b_n(U_{11})$ note that for any $\alpha\in B_n(U_{11})$ we can may write
    $\alpha=(\beta,\gamma)$ where $\gamma=(n,\dots,n-j+1)$ and $\beta$ is any
    permutation of $\aa,\ab,\dots,\overline{n-j}$, so $b_n(U_{11})=0!+1!+\dots+n!$
    for all $n\geq 0$, as claimed in $(4.6)$.

To verify $(4.7)$ note that $\alpha\in B_n(U_2)$ has two cases either
    $\alpha$ contain one symbol unbarred or all the symbols of $\alpha$ are barred. In the
    second case there are $n!$ signed permutations. In the first case, let
    $\alpha=(\beta,i,\gamma)$ where $\beta,\gamma$ permutations of subsets
    of $\aa,\ab,\dots,\an$, but by definition we get that $|\beta_i|>i>|\gamma_j|$, hence
    there are $\sum_{i=1}^n (n-i)!(i-1)!$ signed permutations, which means $(4.7)$ holds.
\end{proof}

\section{More than four signed patterns}
By taking advantage of Proposition \ref{sym}, the question of
determining the values $b_n(T)$ for $n\geq3$ where $T\subset B_2$
and $|T|=5$, for the $56$ choices of five $2$-letters signed
patterns reduces to $10$ cases, as follows: {\small
$$\begin{array}{lll}
W_1=\{12,1\ab,\aa2,\aa\ab,21\};&\ \     W_2=\{12,1\ab,\aa2,\aa\ab,2\aa\};\ \    & W_3=\{12,1\ab,\aa2,21,2\aa\};\\
W_4=\{12,1\ab,\aa2,21,\ab\aa\};&\ \     W_5=\{12,1\ab,\aa2,2\aa,\ab1\};\ \  & W_6=\{12,1\ab,\aa2,2\aa,\ab\aa\};\\
W_7=\{12,1\ab,\aa\ab,21,2\aa\};&\ \     W_8=\{12,1\ab,\aa\ab,21,\ab1\};\ \  & W_9=\{12,1\ab,\aa\ab,21,\ab\aa\};\\
W_{10}=\{12,1\ab,\aa\ab,2\aa,\ab1\}.&                       &
\end{array}$$}

With the aid of a computer we have calculated the cardinality of
$B_n(W_j)$ for $j=1,2,\cdots,10$. From these results we arrived at
the plausible conjecture of Theorem \ref{th5} (some of which are
trivially true).

\begin{theorem}
\label{th5} Given $n$ and set $T$ of $2$-letters signed patterns
such that $|T|=5$. The value $b_n(T)$ satisfies one of the
following relations, according to which orbit (under reversal,
barring, complementation) contains $T$:
\begin{enumerate}
\item   $b_n(W_9)=0$;
\item   $b_n(W_4)=3$;
\item   $b_n(W_1)=b_n(W_2)=b_n(W_6)=b_n(W_7)=b_n(W_8)=b_n(W_{10})=n+1$;
\item   $b_n(W_5)=1+n!$;
\item   $b_n(W_3)=(n+1)(n-1)!$.
\end{enumerate}
\end{theorem}

Additionally, the question of determining
the values $b_n(T)$ where $T\subset B_2$ and $|T|=6$, for the $28$ choices of six
$2$-letters signed patterns reduces to $8$ cases, as follows:
{\small $$\begin{array}{lll}
V_1=\{12,1\ab,\aa2,\aa\ab,21,2\aa\};\, &   V_2=\{12,1\ab,\aa2,\aa\ab,21,\ab\aa\};\, & V_3=\{12,1\ab,\aa2,\aa\ab,2\aa,\ab1)\};\\
V_4=\{12,1\ab,\aa2,21,2\aa,\ab1\};      &   V_5=\{12,1\ab,\aa2,21,2\aa,\ab\aa\};      & V_6=\{12,1\ab,\aa2,2\aa,\ab1,\ab\aa\};\\
V_7=\{12,1\ab,\aa\ab,21,2\aa,\ab\aa\};  &   V_8=\{12,1\ab,\aa\ab,21,\ab1,\ab\aa\}.  &
\end{array}$$}

With the aid of a computer we have calculated the cardinality of
$B_n(V_j)$ for $j=1,2,\cdots,8$. From these results we arrived at
the plausible conjecture of Theorem \ref{th6} (some of which are
trivially true).

\begin{theorem}
\label{th6} Given $n$ and set $T$ of $2$-letters signed patterns
such that $|T|=6$. The value $b_n(T)$ for $n\geq2$ satisfies one
of the following relations, according to which orbit (under
reversal, barring, complementation) contains $T$:
\begin{enumerate}
\item   $b_n(V_2)=b_n(V_7)=b_n(V_8)=0$;
\item   $b_n(V_1)=b_n(V_3)=b_n(V_5)=b_n(V_6)=2$;
\item   $b_n(V_4)=n!$ for all $n\geq 2$.
\end{enumerate}
\end{theorem}

Additionally, the question of determining the values $b_n(T)$
where $T\subset B_2$ and $|T|\geq 7$, for the $9$ choices of seven
or eight $2$-letters signed patterns reduces to $3$ cases which
are:
$$U_1=B_2,\ U_2=\{12,1\ab,\aa2,\aa\ab,21,2\aa,\ab\aa\},\
U_3=\{12,1\ab,\aa2,\aa\ab,21,2\aa,\ab1\}.$$

With the aid of a computer we have calculated the cardinality of
$B_n(U_j)$ for $j=1,2,3$. From these results we arrived at the
plausible conjecture of Theorem \ref{th7} (some of which are
trivially true).

\begin{theorem}\label{th7}
Given $n$ and set $T$ of $2$-letters signed patterns such that
$|T|=7,8$. The value $b_n(T)$ for $n\geq3$ satisfies one of the
following relations, according to which orbit (under reversal,
barring, complementation) contains $T$:
\begin{enumerate}
\item   $b_n(U_1)=b_n(U_2)=0$;

\item   $b_n(U_3)=1$.
\end{enumerate}
\end{theorem}

\end{document}